\documentclass{article}

\usepackage{latexsym}

\usepackage{amsmath, amssymb}

\title{The graph zeta functions with respect to the group matrix of a finite group}
\author{Tsuyoshi MIEZAKI \\
Faculty of Science and Engineering, Waseda University, \\
Tokyo 169-8555, JAPAN \\
email: miezaki@waseda.ac.jp \\
Iwao SATO  \\ 
Oyama National College of Technology, \\ 
Oyama, Tochigi 323-0806, JAPAN \\
e-mail: isato@oyama-ct.ac.jp}

\begin{document}
 \maketitle

\begin{abstract}
In this paper, we present formulas for the edge zeta function and the second weighted zeta 
function with respect to the group matrix of a finite abelian group $\Gamma $. 
Furthermore, we give another proof of Dedekind Theorem for the group determinant of $\Gamma $ 
by the decomposition formula for a matrix of a group covering of a digraph. 
Finally, we treat the weighted complexity of the complete graph with entries of the group matrix 
of $\Gamma $ as arc weights. 
\end{abstract}


\section{Introduction}

Zeta functions of graphs started from the Ihara zeta functions of regular graphs 
by Ihara [5]. 
In [5], Ihara showed that their reciprocals are explicit polynomials. 
A zeta function of a regular graph $G$ associated with a unitary 
representation of the fundamental group of $G$ was developed by 
Sunada [11,12]. 
Hashimoto [3] generalized Ihara's result on the Ihara zeta function of 
a regular graph to an irregular graph, and showed that its reciprocal is 
again a polynomial given by a determinant containing the edge matrix. 
Bass [1] presented another determinant expression for the Ihara zeta function 
of an irregular graph by using its adjacency matrix. 

Weighted zeta functions of a graph are generalizations of the Ihara zeta function of a graph. 
Stark and Terras [10] introduced the edge zeta function of a graph, and gave 
its determinant expression by using the edge matrix. 
Watanabe and Fukumizu [13] presented a determinant expression for the edge zeta function 
of a graph with $n$ vertices by using two $n \times n$ matrices. 
Sato [8] defined the second weighted zeta function of a graph, and presented 
its determinant expression. 

In this paper, we present formulas for the edge zeta function and the second weighted zeta 
function with respect to the group matrix of a finite abelian group.  

This paper is organized as follows: 
In Section 2, we give a short review for the Ihara zeta function, the edge zeta function and 
the second weighted zeta function of a graph. 
Furthermore, we state Dedekind Theorem on the determinant of the group matrix of a finite abelian group.  
In Section 3, we present a formula for the edge zeta function with respect to the group matrix 
of a finite abelian group $\Gamma $. 
In Section 4, we give a formula for the second weighted zeta function with respect to the group matrix 
of $\Gamma $. 
In Section 5, we present another proof of Dedekind Theorem for the group determinant of $\Gamma $ 
by the decomposition formula for a matrix of a group covering of a digraph. 
In Section 6, we treat the weighted complexity of the complete graph with entries of the group matrix 
of $\Gamma $ as arc weights. 

\section{Preliminaries} 

\subsection{Zeta functions of a graph} 

Graphs treated here are finite.
Let $G=(V(G),E(G))$ be a connected graph (possibly multiple edges and loops) 
with the set $V(G)$ of vertices and the set $E(G)$ of unoriented edges $uv$ 
joining two vertices $u$ and $v$. 
For $uv \in E(G)$, an arc $(u,v)$ is the oriented edge from $u$ to $v$. 
Set $D(G)= \{ (u,v),(v,u) \mid uv \in E(G) \} $. 
For $e=(u,v) \in D(G)$, set $u=o(e)$ and $v=t(e)$. 
Furthermore, let $e^{-1}=(v,u)$ be the {\em inverse} of $e=(u,v)$. 

A {\em path $P$ of length $n$} in $G$ is a sequence 
$P=(e_1, \ldots ,e_n )$ of $n$ arcs such that $e_i \in D(G)$,
$t( e_i )=o( e_{i+1} ) \ (1 \leq i \leq n-1)$, 
where indices are treated $mod \  n$. 
 If $e_i =(v_i, v_{i+1} )\  (1 \leq i \leq n)$, then we write 
$P=(v_1, \ldots ,v_{n+1} )$. 
Set $ \mid P \mid =n$, $o(P)=o( e_1 )$ and $t(P)=t( e_n )$. 
Also, $P$ is called an {\em $(o(P),t(P))$-path}. 
We say that a path $P=(e_1, \ldots ,e_n )$ has a {\em backtracking} 
if $ e^{-1}_{i+1} =e_i $ for some $i \ (1 \leq i \leq n-1)$. 
A $(v, w)$-path is called a {\em $v$-cycle} 
(or {\em $v$-closed path}) if $v=w$. 
The {\em inverse cycle} of a cycle 
$C=( e_1, \ldots ,e_n )$ is the cycle 
$C^{-1} =( e^{-1}_n , \ldots ,e^{-1}_1 )$.

We introduce an equivalence relation between cycles. 
Two cycles $C_1 =(e_1, \ldots ,e_m )$ and 
$C_2 =(f_1, \ldots ,f_m )$ are called {\em equivalent} if there exists 
$k$ such that $f_j =e_{j+k} $ for all $j$. 
The inverse cycle of $C$ is in general not equivalent to $C$. 
Let $[C]$ be the equivalence class which contains a cycle $C$. 
Let $B^r$ be the cycle obtained by going $r$ times around a cycle $B$. 
Such a cycle is called a {\em power} of $B$. 
A cycle $C$ is {\em reduced} if both $C$ and $C^2 $ have no backtracking. 
Furthermore, a cycle $C$ is {\em prime} if it is not a power of 
a strictly smaller cycle. 
Note that each equivalence class of prime, reduced cycles of a graph $G$ 
corresponds to a unique conjugacy class of 
the fundamental group $ \pi {}_1 (G,v)$ of $G$ at a vertex $v$ of $G$. 

The {\em Ihara zeta function} of a graph $G$ is 
a function of $u \in \mathbb{C}$ with $|u|$ sufficiently small, defined by 
\[
{\bf Z} (G, u)= {\bf Z}_G (u)= \prod_{[C]} (1- u^{ \mid C \mid } )^{-1} ,
\]
where $[C]$ runs over all equivalence classes of prime, reduced cycles 
of $G$(see [5]).

Let $m$ be the number of edges of $G$. 
Furthermore, let two $2m \times 2m$ matrices 
${\bf B} =( {\bf B}_{ef} )_{e,f \in D(G)} $ and 
${\bf J}_0 =( {\bf J}_{ef} )_{e,f \in D(G)} $ be defined as follows: 
\[
{\bf B}_{ef} =\left\{
\begin{array}{ll}
1 & \mbox{if $t(e)=o(f)$, } \\
0 & \mbox{otherwise, }
\end{array}
\right.
\  
{\bf J}_{ef} =\left\{
\begin{array}{ll}
1 & \mbox{if $f= e^{-1} $, } \\
0 & \mbox{otherwise.}
\end{array}
\right.
\]
Then ${\bf B} -{\bf J}_0 $ is called the {\em edge matrix} of $G$.

\newtheorem{theorem}{Theorem}
\begin{theorem}[Hashimoto; Bass]
Let $G$ be a connected simple graph with $n$ vertices and $m$ edges. 
Then the reciprocal of the Ihara zeta function of $G$ is given by 
\[
{\bf Z} (G, u)^{-1} = \det ( {\bf I}_{2m} - u ( {\bf B} - {\bf J}_0 ))
=(1- u^2 )^{m-n} \det ( {\bf I}_n -u {\bf A} (G)+ u^2 ({\bf D} -{\bf I}_n )), 
\]
where ${\bf A} (G)$ is the adjacency matrix of $G$, and ${\bf D} =( d_{ij} )$ 
is the diagonal matrix with $d_{ii} = \deg v_i $ where $V(G)= \{ v_1 , \ldots , v_n \} $. 
\end{theorem}

Now, we state a weighted zeta function of a graph $G$. 
Stark and Terras [10] defined the edge zeta function of a graph $G$ with $n$ vertices. 
Let $G$ be a connected simple graph and $D(G)= \{ e_1, \ldots, e_m, e_{m+1} , \ldots , 
e_{2m} \} (e_{m+i} = e^{-1}_i (1 \leq i \leq m ))$. 
We introduce $2m$ variables ${\bf u} =( u_1 , \ldots , u_{2m} )$, and set $g(C)=  u_{i_1 } \cdots u_{i_k} $ 
for each cycle $C=( e_{i_1 }, \ldots, e_{i_k} )$ of $G$. 
Then the {\em edge zeta function} $\zeta {}_G (u)$ of $G$ is defined by 
\[
\zeta {}_G ( {\bf u} )= \prod_{[C]} (1-g(C) )^{-1} , 
\]
where $[C]$ runs over all equivalence classes of prime, reduced cycles 
of $G$.

\begin{theorem}[Stark and Terras]
Let $G$ be a connected simple graph with $m$ edges. 
Then 
\[
\zeta {}_G ( {\bf u} ) {}^{-1} = 
\det ( {\bf I}_{2m} - ( {\bf B} - {\bf J} {}_0) {\bf U} ) 
= \det ( {\bf I}_{2m} - {\bf U} ( {\bf B} - {\bf J} {}_0)) ,  
\]
where two matrices ${\bf B} =( {\bf B}_{e,f} )_{e,f \in D(G)} $ and 
${\bf J}_{0} =( {\bf J}_{e,f} )_{e,f \in D(G)} $ are given by  
\[
{\bf B}_{e,f} =\left\{
\begin{array}{ll}
1 & \mbox{if $t(e)=o(f)$, } \\
0 & \mbox{otherwise, }
\end{array}
\right.
\  
{\bf J}_{e,f} =\left\{
\begin{array}{ll}
1 & \mbox{if $f= e^{-1} $, } \\
0 & \mbox{otherwise.}
\end{array}
\right.
\]
and 
\[
{\bf U} = {\rm diag} ( u_1 , \ldots, u_m , u_{m+1} , \ldots , u_{2m} ) . 
\]
\end{theorem}

Watanabe and Fukumizu [13] presented a determinant expression for the edge zeta function 
of a graph $G$ with $n$ vertices by $n \times n$ matrices. 
Then we define an $n \times n$ matrix $\widehat{{\bf A}} = \widehat{{\bf A}} (G)=(a_{xy} )$ as follows: 
\[
a_{xy} =\left\{
\begin{array}{ll}
u_{(x,y)}/(1- u_{(x,y)} u_{(y,x)} ) & \mbox{if$(x,y) \in D(G)$, } \\
0 & \mbox{otherwise. }
\end{array}
\right. 
\]
Furthermore, an $n \times n$ matrix $\widehat{{\bf D}} =\widehat{{\bf D}} (G)=(d_{xy} )$ is the diagonal matrix 
defined by 
\[
d_{xx} = \sum_{o(e)=x} 
\frac{u_e u_{e^{-1}} }{1- u_e u_{ e^{-1}} } . 
\]

\begin{theorem}[Watanabe and Fukumizu]
\[
\zeta {}_G ( {\bf u} ) {}^{-1} = \det ( {\bf I}_{n} + \widehat{{\bf D}} 
- \widehat{{\bf A}} ) \prod^{m}_{i=1} ( 1- u_{e_i} u_{e^{-1}_i } ) . 
\]
\end{theorem}

If $G$ is a multigraph, in Theorem 3, the matrix $\widehat{{\bf A}} =(a_{xy} )$ is given as follows: 
\[
a_{xy} =\left\{
\begin{array}{ll}
\sum \{ u_{(x,y)}/(1- u_{(x,y)} u_{(y,x)} ) \mid (u,v) \in D(G) \} & \mbox{if$(x,y) \in D(G)$, } \\
0 & \mbox{otherwise. }
\end{array}
\right. 
\]

Next, we state a weighted zeta function of a graph $G$ different from the edge zeta function. 
Let $G$ be connected simple graph with $n$ vertices and $m$ edges. 
Then an $n \times n$ matrix ${\bf W} (G)=( w_{uv} )$ is given as follows: 
\[
w_{uv} =\left\{
\begin{array}{ll}
nonzero \ complex \ number & \mbox{if $(u,v) \in D(G)$, } \\
0 & \mbox{otherwise. }
\end{array}
\right. 
\] 
The matrix ${\bf W} (G)$ is the {\em weighted matrix} of $G$. 
Set $w(u,v)=w_{uv} , \ u,v \in V(G)$ and $w(e)= w_{uv} , \ e=(u,v) \in D(G)$. 

Furthermore, we define a $2m \times 2m$ matrix ${\bf B}_w =( B^{(w)}_{ef} )_{e,f \in D(G)} $ 
be defined as follows: 
\[
B^{(w)}_{ef} =\left\{
\begin{array}{ll}
w(f) & \mbox{if $t(e)=o(f)$, } \\
0 & \mbox{otherwise, }
\end{array}
\right.
\]
Then the {\em second weighted zeta function} of a graph $G$ is defined as follows([8]): 
\[
{\bf Z}_1 (G,w,u)= \det ( {\bf I}_{2m} -u( {\bf B}_w - {\bf J}_0 ))^{-1} . 
\]

If $w= {\bf 1}$, i.e., $w(e)=1 $ for each $e \in D(G)$, then the second weighted zeta function is equal to 
the Ihara zeta function: 
\[
{\bf Z}_1 (G,w,u)= {\bf Z} (G,u) . 
\]

The determinant expression of Ihara type for the second weighted zeta function is given as follows([8]):

\begin{theorem}[Sato]
Let $G$ be connected simple graph with $n$ vertices and $m$ edges, and ${\bf W} (G)$ a weighted matrix of $G$. 
Then the reciprocal of the second weighted zeta function of $G$ is 
\[
{\bf Z}_1 (G,w,u )^{-1} =(1- u^2 )^{m-n} 
\det ({\bf I}_n -u {\bf W} (G)+ u^2 ( {\bf D}_w - {\bf I}_n )) , 
\]
where the matrix ${\bf D}_w =( d_{uv} )$ is an $n \times n$ diagonal matrix with 
\[
d_{uu} = \sum_{o(e)=u} w(e) . 
\]
\end{theorem}

If $G$ is a multigraph, in Theorem 4, the matrix ${\bf W} =(w_{xy} )$ is given as follows: 
\[
w_{xy} =\left\{
\begin{array}{ll}
\sum \{ w(u,v) \mid (u,v) \in D(G) \} & \mbox{if$(x,y) \in D(G)$, } \\
0 & \mbox{otherwise. }
\end{array}
\right. 
\]

\subsection{The group matrix} 

Let $\Gamma = \{ g_1 =1, g_2 , \ldots , g_n \} $ be a finite group and $R$ a polynomial ring containing 
the indeterminates $x_g $ for $g \in \Gamma $. 
Then the {\em group matrix} ${\bf M} ( \Gamma )$ of $\Gamma $ is defined as follows: 
\[
({\bf M} ( \Gamma ))_{g_i g_j } = x_{g^{-1}_i g_j } \ (1 \leq i,j \leq n) . 
\]

The determinant of the group matrix ${\bf M} ( \Gamma )$ for an abelian group $ \Gamma $ is presented as follows:

\begin{theorem}[Dedekind] 
Let $ \Gamma $ be a finite abelian group. 
Then 
\[
\det ({\bf M} ( \Gamma ))= \prod_{\chi \in \hat{\Gamma } } ( \sum_{g \in \Gamma } \chi (g) x_g ) , 
\] 
where $\chi $ runs over all inequivalent characters of $\Gamma $. 
\end{theorem}

\section{The edge zeta function with respect to the group matrix} 

Let $\Gamma = \{ g_1 =1, g_2 , \ldots , g_n \} $ be a finite abelian group, $R$ a polynomial ring containing 
the indeterminates $x_g $ for $g \in \Gamma $ and ${\bf M} ( \Gamma )$ the group matrix of $\Gamma $. 
Then we define the complete graph $K_{\Gamma } $ with $n$ loops as follows: 
\[
V(K_{\Gamma } )= \Gamma \ and \ E( K_{\Gamma } )= \{ gh \mid g,h \in \Gamma \} . 
\]
Furthermore, the arc set of the symmetric digraph $D_{K_{\Gamma }} $ is 
\[
D( K_{\Gamma } )= \{ (g,h) \mid g.h \in \Gamma \} . 
\]
We introduce an arc weight $u: D(K_{\Gamma } ) \longrightarrow R$ as follows: 
\[
u_e = x_{g^{-1} h} \ for \ e=(g,h) \in D(K_{\Gamma } ) . 
\]

Then the following result follows.

\begin{theorem} 
Let $\Gamma = \{ g_1 =1, g_2 , \ldots , g_n \} $ be a finite abelian group, $R$ a polynomial ring containing 
the indeterminates $x_g $ for $g \in \Gamma $ and $K_{\Gamma } $ the complete graph with $n$ loops for $\Gamma $. 
Furthermore, let 
\[
y_1 =1+ \frac{2 x^2_1 }{1- x^2_1 } + \sum^n_{j=2} \frac{x_{g_j} x_{g^{-1}_j }}{1- x_{g_j} x_{g^{-1}_j }} , \ 
y_{g^{-1}_i g_j } =- \frac{(1+ \delta {}_{ij} ) x_{g^{-1}_i g_j} }{1- x_{g^{-1}_i g_j} x_{g_i g^{-1}_j }} . 
\]
Then 
\[
\zeta {}_{K_{\Gamma }} ({\bf u} )^{-1} = \prod_{g \in \Gamma } (1- x_g x_{g^{-1}} )^n \prod_{\chi } ( \sum_{g \in \Gamma } \chi (g) y_g ) . 
\]
\end{theorem}

{\bf Proof}.  Let $m=n(n-1)/2+n=n(n+1)/2$ and $E( K_{\Gamma } )= \{ f_1 , \ldots, f_n , f_{n+1} , \ldots , f_m \}$ such that 
\[
f_i= g_i \ (1 \leq i \leq n) ; \ f_j = g_{j_1 } g_{j_2 } \ (j=n+, 1 \ldots , m; 1 \leq j_1 < j_2 \leq n) . 
\]
Furthermore, let $D( K_{\Gamma } )= \{ e_1 , \ldots , e_n , e_{n+1} , \ldots , e_m , 
e^{-1}_1 , \ldots , e^{-1}_n , e^{-1}_{n+1} , \ldots , e^{-1}_m \} $ such that 
\[
e_i =(g_i , g_i ) \ (1 \leq i \leq n) ; \ e_j = ( g_{j_1 } , g_{j_2 } ) \ (j=n+, 1 \ldots , m; 1 \leq j_1 < j_2 \leq n) . 
\] 
Moreover, let $u_{e_j } =x_{g^{-1}_{j_1} g_{j_2 } } $ for $e_j =(g_{j_1 } , g_{j_2 } ) \in D( K_{\Gamma } )$. 
Then, by Theorem 3, we have 
\[
\zeta {}_{K_{\Gamma }} ({\bf u} )^{-1} = \prod^m_{j=1} (1- u_{e_j } u_{e^{-1}_j } ) \det ({\bf I}_n + \widehat{{\bf D}} - \widehat{{\bf A}} ) . 
\]
Here, two matrices $\widehat{{\bf A}} =( A_{g_i g_j } )$ and $\widehat{{\bf D}} =( D_{g_i g_j } )$ are given as follows: 
\[
A_{g_i g_j } = \frac{(1+ \delta {}_{ij} ) x_{g^{-1}_i} x_{g_j }}{1- x_{g^{-1}_i g_j } x_{g_i g^{-1}_j }} , 
\ 
D_{g_i g_i } = \sum^m_{j=1} \frac{(1+ \delta {}_{ij} ) x_{g^{-1}_i g_j } x_{g_i g^{-1}_j }}{1- x_{g^{-1}_i g_j } x_{g_i g^{-1}_j } } . 
\]

But, we have 
\[
\{ g^{-1}_i g_j \mid 1 \leq j \leq n \} = \{ g_1 , \ldots , g_n \} = \Gamma . 
\]
Then we get  
\[ 
D_{g_i g_i } = \frac{2 x^2_1 }{1- x^2_1 } + \sum_{g \in \Gamma , g \neq 1} \frac{x_g x_{g^{-1}} }{1- x_g x_{g^{-1} } } = {\rm const} . 
\]
Thus, let 
\[
y_1 = y_{g_1 } =1+ D_{g_i g_i } ; \ y_{g^{-1}_i g_j } =- A_{g_i g_j } . 
\]
Furthermore, let the $n \times n$ matrix ${\bf N} =( N_{ij} )$ be defined as follows: 
\[
N_{ij} = y_{g^{-1}_i g_j } . 
\]
Then the matrix ${\bf N} $ is a new group matrix of $\Gamma $. 
Therefore, it follows that 
\[
\zeta {}_{K_{\Gamma }} ({\bf u} )^{-1} = \prod^m_{j=1} (1- u_{e_j } u_{e^{-1}_j } ) \det ({\bf N} ) . 
\] 
By Theorem 5, we have 
\[
\zeta {}_{K_{\Gamma }} ({\bf u} )^{-1} = \prod_{g \in \Gamma } (1- x_g x_{g^{-1}} )^n \prod_{\chi } ( \sum_{g \in \Gamma } \chi (g) y_g ) . 
\]
$\Box$

\section{The second weighted zeta function with respect to the group matrix} 

Let $\Gamma = \{ g_1 =1, g_2 , \ldots , g_n \} $ be a finite abelian group, $R$ a polynomial ring containing 
the indeterminates $x_g $ for $g \in \Gamma $ and ${\bf M} ( \Gamma )$ the group matrix of $\Gamma $. 
Furthermore, let $K_{\Gamma } $ be the complete graph with $n$ loops with respect to $\Gamma $. 

Then the following result follows.

\begin{theorem} 
Let $\Gamma = \{ g_1 =1, g_2 , \ldots , g_n \} $ be an finite abelian group, $R$ a polynomial ring containing 
the indeterminates $x_g $ for $g \in \Gamma $, ${\bf M} ( \Gamma )$ the group matrix of $\Gamma $ and 
$K_{\Gamma } $ the complete graph with $n$ loops for $\Gamma $. 
Furthermore, let the weighted matrix ${\bf W} ( K_{\Gamma } )$ be given as follows: 
\[
{\bf W} ( K_{\Gamma } )= {\bf M} ( \Gamma ) . 
\] 
Then 
\[
{\bf Z}_1 ( K_{\Gamma } ,w,u )^{-1} =(1- u^2 )^{n(n-1)/2} \prod_{\chi } (1-2 x_1 u-u \sum_{g \neq 1} \chi (g) x_g + u^2 ( \sum_{g \in \Gamma } x_g +x_1 -1)) . 
\]
\end{theorem}

{\bf Proof}.  Since ${\bf W} ( K_{\Gamma )} = {\bf M} ( \Gamma )$, the weight function $w: D( K_{\Gamma } ) \longrightarrow R$ is 
given as follows: 
\[
w(e)= x_{g^{-1}_i g_j } , \ e=(g_i , g_j ) .
\]
Furthermore, we have $m=|E( K_{\Gamma } )|=n(n-1)/2+n=n(n+1)/2$. 
By Theorem 4, we have 
\[
\begin{array}{rcl} 
{\bf Z}_1 ( K_{\Gamma } ,w,u)^{-1} & = & (1- u^2 )^{n(n+1)/2 -n} \det ( {\bf I}_n -u {\bf M} ( {\Gamma } )+u^2 ( {\bf D}_w -{\bf I}_n )) \\ 
\  &   &                \\ 
\  & = & (1- u^2 )^{n(n-1)/2} \det ( {\bf I}_n -u {\bf M} ( {\Gamma } )+u^2 ( {\bf D}_w -{\bf I}_n )) . 
\end{array} 
\]

But, we have 
\[
( {\bf D}_w )_{g_i g_i } = \sum_{o(e)= g_i } w(e)= \sum^n_{j=1} (1+ \delta {}_{ij} ) x_{g^{-1}_i g_j } 
= \sum_{g \in \Gamma } x_g + x_1 = \mu \ ({\rm const} ) . 
\]
Thus, let 
\[
y_1 =1-2 x_1 u+ u^2 ( \mu -1)=1-2 x_1 u + u^2 ( \sum_{g \in \Gamma } x_g + x_1 -1 ); \ 
y_g = - x_g u . 
\]
Furthermore, let the $n \times n$ matrix ${\bf Y} =( Y_{ij} )$ be defined as follows: 
\[
{\bf Y} =( y_{g^{-1}_i g_j } )_{1 \leq i,j \leq n} . 
\]
Moreover, the matrix ${\bf N} $ is a new group matrix of $\Gamma $. 
Then we have  
\[
{\bf Z}_1 ( K_{\Gamma } ,w,u)^{-1} =(1- u^2 )^{n(n-1)/2} \det ( {\bf Y} ) .  
\] 
By Theorem 5, we have 
\[ 
\begin{array}{rcl} 
{\bf Z}_1 ( K_{\Gamma } ,w,u)^{-1} & = & (1- u^2 )^{n(n-1)/2} \prod_{\chi } ( \sum_{g \in \Gamma } \chi (g) y_g ) \\
\  &   &                \\ 
\  & = & (1- u^2 )^{n(n-1)/2} \prod_{\chi } ( \chi (1) y_1 + \sum_{g \neq 1} \chi (g) y_g ) \\
\  &   &                \\ 
\  & = & (1- u^2 )^{n(n-1)/2} \prod_{\chi } (1-2 x_1 u+ u^2 ( \sum_g x_g + x_1 -1)-u \sum_{g \neq 1} \chi (g) x_g ) . 
\end{array} 
\] 
$\Box$

\section{Another proof for Dedekind Theorem} 

At first, we consider a group covering of a digraph. 

Let $D$ be a digraph and $ \Gamma $ a finite group.
Then a mapping $ \alpha : A(D) \longrightarrow \Gamma $
is called an {\em pseudo ordinary voltage} {\em assignment}
if $ \alpha (v,u)= \alpha (u,v)^{-1} $ for each $(u,v) \in A(D)$ 
such that $(v,u) \in A(D)$.
The pair $(D, \alpha )$ is called an 
{\em ordinary voltage digraph}.
The {\em derived digraph} $D^{ \alpha } $ of the ordinary
voltage digraph $(D, \alpha )$ is defined as follows:
$V(D^{ \alpha } )=V(D) \times \Gamma $ and $((u,h),(v,k)) \in 
A(D^{ \alpha })$ if and only if $(u,v) \in A(D)$ and $k=h \alpha (u,v) $. 
The {\em natural projection} 
$ \pi : D^{ \alpha } \longrightarrow D$ is defined by 
$ \pi (u,h)=u$. 
The graph $D^{ \alpha }$ is called a 
{\em derived digraph covering} of $D$ with voltages in 
$ \Gamma $ or a {\em $ \Gamma $-covering} of $D$.
A $\Gamma $-covering of the symmetric digraph of a graph $G$ is a regular covering of $G$ (see [2]). 

In the $\Gamma $-covering $D^{ \alpha } $, set $v_g =(v,g)$ and $e_g =(e,g)$, 
where $v \in V(D), e \in A(D), g\in \Gamma $. 
For $e=(u,v) \in A(D)$, the arc $e_g$ emanates from $u_g$ and 
terminates at $v_{g \alpha (e)}$. 
Note that $ e^{-1}_g =(e^{-1} )_{g \alpha (e)}$. 

Let $D$ be a connected digraph with $n$ vertices and $m$ arcs, 
$ \Gamma $ a finite group, $ \alpha : A(D) \longrightarrow \Gamma  $ 
a pseudo ordinary voltage assignment and ${\bf W} = {\bf W} (D)$ a weighted matrix of $D$. 
Then we define the {\em weighted matrix} $\tilde{{\bf W}} = {\bf W} ( D^{\alpha } )=( \tilde{w} ( u_g , v_h ))$ 
{\em derived from} ${\bf W} $ as follows: 
\[
\tilde{w} ( u_g , v_h ) :=\left\{
\begin{array}{ll}
w(u,v) & \mbox{if $(u,v) \in A(D)$ and $h=g \alpha (u,v)$, } \\
0 & \mbox{otherwise.}
\end{array}
\right.
\]

Now, we consider the determinant $\det ( {\bf W} (D^{\alpha } ))$ of the weighted matrix ${\bf W} (D^{\alpha } )$ of 
a group covering $D^{\alpha } $ of $G$. 

For $g \in \Gamma $, let $n \times n$ matrices ${\bf W}_{g} ={\bf W}_{g} (D)=( W^{(g)}_{uv} )_{u,v \in V(D)} $ is 
defined as follows: 
\[  
W^{(g)}_{uv} =\left\{
\begin{array}{ll}
w(u,v) & \mbox{if $(u,v) \in A(D)$ and $\alpha (u,v)=g$, } \\
0 & \mbox{otherwise. }
\end{array}
\right.
\]

Let \( {\bf M}_{1} \oplus \cdots \oplus {\bf M}_{s} \) be the 
block diagonal sum of square matrices 
${\bf M}_{1}, \cdots , {\bf M}_{s}$. 
If \( {\bf M}_{1} = {\bf M}_{2} = \cdots = {\bf M}_{s} = {\bf M} \),
then we write 
\( s \circ {\bf M} = {\bf M}_{1} \oplus \cdots \oplus {\bf M}_{s} \).

\begin{theorem} 
Let $D$ be a connected digraph with $n$ vertices and $m$ arcs, 
$ \Gamma $ a finite group, $ \alpha : A(D) \longrightarrow \Gamma  $ 
a pseudo ordinary voltage assignment, and ${\bf W} $ a weighted matrix of $D$. 
Set $\mid \Gamma \mid =p$. 
Furthermore, let $ {\rho}_{1} =1, {\rho}_{2} , \ldots , {\rho}_{k} $
be the irreducible representations of $ \Gamma $, and 
$d_i$ the degree of $ {\rho}_{i} $ for each $i$, where 
$d_1=1$.
Suppose that the $ \Gamma $-covering $D^{ \alpha } $ of $D$ is connected. 
Then the determinant $\det ( {\bf W} (D^{\alpha } ))$ of the weighted matrix of $G^{\alpha } $ is
\[
\det ( {\bf W} (D^{\alpha } ))= \prod^k_{j=1} \det ( \sum_{h \in \Gamma } \rho {}_j (h) \bigotimes {\bf W}_h )^{d_j} . 
\]  
\end{theorem}

{\bf Proof }.  Let $V(D)= \{ v_1, \ldots , v_{n} \} $ and 
$ \Gamma = \{ 1=g_1, g_2, \ldots ,g_p \} $.
Arrange vertices of $D^{ \alpha } $ in $n$ blocks:
$(v_1,1), \ldots , (v_{n},1);(v_1,g_2), \ldots , (v_{n},g_2); 
\cdots ; (v_1,g_p), \ldots ,(v_{n},g_p). $
We consider the matrix ${\bf W} ( D {}^{ \alpha } ) $ under this order.

For $h \in \Gamma $, the matrix ${\bf P}_{h}=(p^{(h)}_{ij} )$ 
is defined as follows:
\[
p^{(h)}_{ij} = \left\{
\begin{array}{ll}
1 & \mbox{if $g_i h=g_j$,} \\
0 & \mbox{otherwise.}
\end{array}
\right.
\]

Next, suppose that $p^{(h)}_{ij} =1 $, i.e., $g_j=g_ih$.
Then $((u,g_i),(v,g_j)) \in A(D {}^{ \alpha } ) $
if and only if $(u,v) \in A(D)$ and 
$g_{j} = g_{i} \alpha (u,v)$,
i.e., $ \alpha (u,v)=g^{-1}_{i} g_j =g^{-1}_{i} g_i h=h$. 
Furthermore, if $((u,g_i),(v,g_j)) \in D(G {}^{ \alpha } ) $, then we have 
\[
( {\bf W} (D^{\alpha } ))_{u_{g_i} v_{g_j }} = \tilde{{\bf w}} (u_{g_i} , v_{g_j} )= {\bf w} (u,v) . 
\]
Thus we have
\[
{\bf W} (D {}^{ \alpha } )= \sum_{h \in \Gamma } {\bf P}_{h} \bigotimes {\bf W} {}_{h} . 
\]

Let $\rho$ be the right regular representation of $ \Gamma $.
Furthermore, let $ {\rho}_{1} =1, {\rho}_{2} , \ldots , {\rho}_{k} $
be all inequivalent irreducible representations of $ \Gamma $, and 
$d_i$ the degree of $ {\rho}_{i} $ for each $i$, where 
$d_1=1$.
Then we have $\rho (h)= {\bf P}_{h} $ for $h \in \Gamma $.
Furthermore, there exists a nonsingular matrix ${\bf P}$ such that
${\bf P}^{-1} \rho (h) {\bf P} = (1) \oplus d_2 \circ {\rho}_{2} (h) 
\oplus \cdots \oplus d_k \circ {\rho}_{k} (h)$ 
for each $h \in \Gamma $(see [9]).  

Putting 
${\bf H} =( {\bf P}^{-1} \bigotimes {\bf I} ) {\bf W} (D {}^{ \alpha } ) 
( {\bf P} \bigotimes {\bf I} )$,
we have 
\[
\begin{array}{rcl} 
{\bf H} & = & \sum_{h \in \Gamma } \{ (1) \oplus d_2 \circ {\rho}_{2} (h) \oplus \cdots \oplus 
d_k \circ {\rho}_{k} (h) \} \bigotimes {\bf W}_{h} \\ 
\  &   &                \\ 
\  & = & \oplus^k_{j=1} ( d_j \circ (\sum_{h \in \Gamma } {\rho}_{j} (h) \bigotimes {\bf W}_{h} )) . 
\end{array} 
\]

Note that ${\bf W} (D) = \sum_{h \in \Gamma } {\bf W} {}_h $ 
and $1+ d^2_2 + \cdots + d^2_k =p$.
Therefore it follows that 
\[
\begin{array}{rcl} 
\det ( {\bf W} (D^{\alpha } )) & = & \prod^k_{i=1} \det ( \sum_{h \in \Gamma } {\rho}_{i} (h) \bigotimes {\bf W} {}_{h} )^{d_i} \\
\  &   &                \\ 
\  & = & \det ( {\bf W} (D)) \prod^k_{i=2} \det ( \sum_{h \in \Gamma } {\rho}_{i} (h) \bigotimes {\bf W} {}_{h} )^{d_i} . 
\end{array} 
\]
$\Box$

In the case that $\Gamma $ is an abelian group, the following result holds.

\newtheorem{corollary}{Corollary}
\begin{corollary}
Let $D$ be a connected digraph, $ \Gamma $ a finite abelian group, 
$ \alpha : A(D) \longrightarrow \Gamma $ a pseudo ordinary voltage 
assignment and ${\bf W} (D)$ a weighted matrix of $D$. 
Then we have 
\[
\det ( {\bf W} (D^{\alpha } ))= \prod_{\chi } \det ( \sum_{h \in \Gamma } \chi (h) {\bf W} {}_{h} ) ,
\]
where $\chi $ runs over all inequivalent characters of $\Gamma $. 
\end{corollary}

Another proof for Dedekind Theorem (Theorem 5): 

Let $\Gamma = \{ g_1 =1, g_2 , \ldots , g_n \} $ be an finite abelian group, $R$ a polynomial ring containing 
the indeterminates $x_g $ for $g \in \Gamma $ and ${\bf M} ( \Gamma )$ the group matrix of $\Gamma $. 
Furthermore, let $D$ be a digraph (bouquet) with one vertex $v$ and $n$ loops $e_1 , \ldots , e_n $. 
Then, let $w: A(D) \longrightarrow R$ be defined as follows: 
\[
w(e_i )= x_{g_i } \ (1 \leq i \leq n) . 
\]

Next, we define a pseudo ordinary assignment $ \alpha : A(D) \longrightarrow \Gamma $ such that 
\[
\alpha (e_i )= g_i \ (1 \leq i \leq n) . 
\]
Then the $\Gamma $-covering $G^{\alpha } $ is as follows: 
\[
V(D^{\alpha } )= \{ (v, g_i ) \mid 1 \leq i \leq n \} \equiv \Gamma , 
\]
and $((v, g_i ), (v, g_j )) \in A(D^{\alpha } )$ if and only if $(v,v) \in A(D)$ and $g_j = g_i \alpha (v,v)$, 
i.e., $(v,v) \in A(D)$ and $\alpha (v,v)= g^{-1}_i g_j $. 
Thus, $D^{\alpha } $ is isomorphic to the complete digraph with $n$ vertices, and $D^{\alpha } $ has only one loop at 
each vertex. 

Furthermore, the weighted matrix $D$ is given as follows: 
\[
{\bf W} (D)=( x_{g_1 } +x_{g_2 } + \cdots + x_{g_n } ) . 
\]
Then the weighted matrix ${\bf W} (D^{\alpha } )=( \tilde{w} ( v_{g_i } , v_{g_j } ))$ of $D^{\alpha } $ is given as follows: 
\[
\tilde{w} ( v_{g_i } , v_{g_j } )= 
\left\{
\begin{array}{ll}
w(v,v) & \mbox{if $(v,v) \in A(D)$ and $g_j = g_i \alpha (v,v)$,} \\
0 & \mbox{otherwise.}
\end{array}
\right.
\] 
Thus, $\alpha (v,v)= g^{-1}_i g_j $ if and only if $w(v,v)= x_{g^{-1}_i g_j } $. 
That is, 
\[
( {\bf W} (D^{\alpha } ))_{g_i g_j } = x_{g^{-1}_i g_j } . 
\] 
Therefore, it follows that 
\[
{\bf W} (D^{\alpha } )= {\bf M} ( \Gamma ) . 
\]
By Corollary 1, we get  
\[
\det ( {\bf W} (D^{\alpha } )= \prod_{\chi } \det ( \sum_{h \in \Gamma } \chi (h) {\bf W}_h ) . 
\]

But, we have 
\[
{\bf W}_h = w(v,v) \ if \ (v,v) \in A(D) \ and \ \alpha (v,v)=h, 
\] 
i.e., 
\[
{\bf W}_h = x_h \ for \ h \in \Gamma .  
\]
Hence, 
\[
\det ({\bf M} ( \Gamma ))= \det ( {\bf W} (D^{\alpha } ))= \prod_{\chi } ( \sum_{h \in \Gamma } \chi (h) x_h ) . 
\]
$\Box$

In the case that $\Gamma $ is not abelian, by Theorem 6 and the above proof for Dedekind Theorem, we obtain the following result.

\begin{corollary} 
Let $\Gamma = \{ g_1 =1, g_2 , \ldots , g_n \} $ be an finite nonabelian group, $R$ a polynomial ring containing 
the indeterminates $x_g $ for $g \in \Gamma $ and ${\bf M} ( \Gamma )$ the group matrix of $\Gamma $. 
Then we have 
\[
\det ( {\bf M} ( \Gamma ))= \prod_{\rho } \det ( \sum_{h \in \Gamma } \rho (h) x_h )^{\deg \rho } ,  
\]
where $\rho $ runs over all inequivalent irreducible representations of $\Gamma $. 
\end{corollary}

\section{The weighted complexity of the complete graph with weights with respect to a group matrix}

The {\em complexity} $ \kappa (G)$($=$ the number of spanning trees in $G$ 
of a connected graph $G$ is closely related to the zeta function of $G$. 
The complexities for various graphs were given in [3,4,7].
Hashimoto expressed the complexity of a regular graph as a limit 
involving its zeta function in [3]. 
For an irregular graph $G$, Hashimoto [4] and Northshield [7] gave 
the value of $(1-u)^{-r} {\bf Z}_G (u)^{-1} $ at $u=1$ in term of 
the complexity of $G$, where $r$ is the Betti number of $G$.

\begin{theorem}[Hashimoto; Northshield]
For any finite graph $G$ such that $n=|V(G)|$, $m=|E(G)|$ and $r=m-n+1>1$, we have
\[
(1-u)^{-r} {\bf Z}_G (u)^{-1} \mid {}_{u=1} = 2^r \chi (G) \kappa (G) ,
\]
where $ \chi (G)=1-r=n-m$ is the Euler number of $G$. 
\end{theorem}

For a connected graph $G$ with $n$ vertices, let 
\[
f_G (u)= \det( {\bf I}_n -u {\bf A} (G)+ u^2 ({\bf D} - {\bf I}_n )) . 
\]

For a graph $G$, Northshield [7] showed that the complexity 
of $G$ is given by the derivative of the function above.

\begin{theorem}[Northshield]
For a connected graph $G$, 
\[
f^{ \prime }_G (1)=2(m-n) \kappa (G ) , 
\]
where $n= \mid V(G) \mid $ and $m= \mid E(G) \mid $. 
\end{theorem}

Let $G$ be a connected graph with $n$ vertices and $m$ edges, and 
${\bf W} = {\bf W} (G)$ a weighted matrix of $G$. 
Then, let 
\[
f_G (w,u)= \det ( {\bf I}_n -u {\bf W} + ( {\bf D}_w - {\bf I}_n ) u^2 ) .
\]

In the case that $w$ is symmetric, i.e., $w(e^{-1} )=w(e)$ 
for each $e \in D(G)$, we consider all spanning arborescences of $G$ 
rooted at any fixed vertex $v \in V(G)$.
A {\em spanning arborescence} $T_v$ of $G$ rooted at $v \in V(G)$ is a 
subdigraph of the symmetric digraph  $D_G$ of $G$ such that ${\rm outdeg} x=1$ for any $x \in V(G)$ 
except $v$, where the outdegree ${\rm outdeg} \ x$ of $x$ is the number of arcs outgoing from $x$ in $D_G$.  
Then, let 
\[
w(T_v)= \prod_{e \in D(T_v)} w(e) . 
\]
Furthermore, let 
\[
\kappa {}_w (G)= \sum_{T_v} w(T_v) . 
\]
The sum of the product of weights of all arcs in those spanning 
arborescences of $G$ are not depended on a vertex $v$ of $G$. 
Then this sum is called the {\em weighted complexity} of $G$. 
Mizuno and Sato [6] showed the following result.

\begin{theorem}[Mizuno and Sato]
\[
\kappa {}_w (G) =\frac{1}{2(w(G)-n)} f^{\prime} (w,1) , 
\]
where $w(G)= \sum_{vw \in E(G)} w(v,w)$. 
\end{theorem}

Let $\Gamma = \{ g_1 =1, g_2 , \ldots , g_n \} $ be an finite abelian group, $R$ a polynomial ring containing 
the indeterminates $x_g $ for $g \in \Gamma $, ${\bf M} ( \Gamma )$ the group matrix of $\Gamma $ and 
$K_{\Gamma } $ the complete graph with $n$ loops for $\Gamma $. 
We consider the weighted complexity of the complete graph $K_{\Gamma }$ for ${\bf W} ( K_{\Gamma })= {\bf M} ( \Gamma )$.

Then the following result follows.

\begin{theorem} 
Let $\Gamma = \{ g_1 =1, g_2 , \ldots , g_n \} $ be an finite abelian group, $R$ a polynomial ring containing 
the indeterminates $x_g $ for $g \in \Gamma $, $K_{\Gamma } $ the complete graph with $n$ loops for $\Gamma $ 
and ${\bf W} ( K_{\Gamma )} = {\bf M} ( \Gamma )$. 
Suppose that $x_g =x_{g^{-1}} , \ g \neq 1 \in \Gamma $ and $x_1 =0 $.  
Then 
\[
\kappa {}_w ( K_{\Gamma } )= \frac{1}{n} \prod_{\chi \neq 1} \sum_{g \neq 1 \in \Gamma } (1- \chi (g)) x_g . 
\]
\end{theorem}

{\bf Proof}.  Since ${\bf W} ( K_{\Gamma )} = {\bf M} ( \Gamma )$, the arc weight $w: D(K_{\Gamma } ) \longrightarrow R$ 
is given as follows: 
\[
w(g_i , g_j )= x_{g^{-1}_i g_j } . 
\]

By Theorems 7, 11, we have 
\[
\kappa {}_w ( K_{\Gamma } ) = \frac{1}{2 \sum_{e \in E(K_{\Gamma } )} w(e)-n} ( \prod_{\chi } 
(1-u \sum_g \chi (g) x_g + u^2 (\sum_g x_g -1))^{\prime } |_{u=1} . 
\]
But, 
\[
\sum_{e \in E(K_{\Gamma } )} w(e)= \frac{n}{2} \sum_{g \neq 1} x_g . 
\]
Furthermore, we have 
\[
\begin{array}{rcl} 
\  &   & ( \prod_{\chi } 
(1-u \sum_g \chi (g) x_g + u^2 (\sum_g x_g -1))^{\prime } \\ 
\  &   &                \\ 
\  & = & \{ (1-u \sum_{g \neq 1} x_g + u^2 (\sum_{g \neq 1} x_g -1))  
\prod_{\chi \neq 1} (1-u \sum_{g \neq 1} \chi (g) x_g + u^2 (\sum_{g \neq 1} x_g -1)) \}^{\prime } \\ 
\  &   &                \\ 
\  & = & \{ (- \sum_{g \neq 1} x_g +2u(\sum_{g \neq 1} x_g -1))  
\prod_{\chi \neq 1} (1-u \sum_{g \neq 1} \chi (g) x_g + u^2 (\sum_{g \neq 1} x_g -1)) \\ 
\  &   &                \\ 
\  & + & (1-u \sum_{g \neq 1} x_g + u^2 (\sum_{g \neq 1} x_g -1))  
( \prod_{\chi \neq 1} (1-u \sum_{g \neq 1} \chi (g) x_g + u^2 (\sum_{g \neq 1} x_g -1)))^{\prime } . 
\end{array} 
\]

Thus, we have 
\[
\begin{array}{rcl} 
\  &   & ( \prod_{\chi } (1-u \sum_g \chi (g) x_g + u^2 (\sum_g x_g -1))^{\prime } |_{u=1} \\  
\  &   &                \\ 
\  & = & \{ (- \sum_{g \neq 1} x_g +2 (\sum_{g \neq 1} x_g -1))  
\prod_{\chi \neq 1} (1- \sum_{g \neq 1} \chi (g) x_g +(\sum_{g \neq 1} x_g -1))+0 \\ 
\  &   &                \\ 
\  & = & ( \sum_{g \neq 1} x_g -2) \prod_{\chi \neq 1} ( \sum_{g \neq 1} (1- \chi (g)) x_g ) . 
\end{array} 
\] 
Therefore, it follows that 
\[
\begin{array}{rcl} 
\kappa {}_w ( K_{\Gamma } ) & = & \frac{1}{2(n/2 \sum_{g \neq 1} x_g -n)} ( \sum_{g \neq 1} x_g -2) 
\prod_{\chi \neq 1} ( \sum_{g \neq 1} (1- \chi (g)) x_g ) \\ 
\  &   &                \\ 
\  & = & \frac{1}{n} \prod_{\chi \neq 1} ( \sum_{g \neq 1} (1- \chi (g)) x_g ) . 
\end{array} 
\] 
$\Box$

In the case that $x_g =1$ for $g \neq 1 \in \Gamma $, we the formula for the complexity of the complete graph.

\begin{corollary}
Let $K_n $ be a the complete graph with $n$ vertices. 
Then 
\[
\kappa (K_n )= n^{n-2} . 
\]
\end{corollary}

{\bf Proof}.  If $x_g =1$ for $g \neq 1 \in \Gamma $, then we have 
\[
\kappa {}_w ( K_{\Gamma } ) = \kappa (K_n ) . 
\]
By Theorem 12, we get 
\[
\begin{array}{rcl} 
\kappa ( K_n ) & = & \frac{1}{n} \prod_{\chi \neq 1} (n-1- \sum_{g \neq 1}  \chi (g)) \\ 
\  &   &                \\ 
\  & = & \frac{1}{n} \prod_{\chi \neq 1} (n- \sum_{g \in \Gamma }  \chi (g)) . 
\end{array} 
\] 
Since $\sum_{g \in \Gamma }  \chi (g)=0$ for $\chi \neq 1$, we obtain 
\[
\kappa ( K_n )= \frac{1}{n} n^{n-1} = n^{n-2} . 
\]
$\Box$

\section*{Statements and Declarations}

\subsection*{Funding}

The second author is supported by JSPS KAKENHI (22K03277).

\subsection*{Competing Interests}

The authors have no affiliation with any organization with a direct or indirect financial interest in the subject matter discussed in the manuscript

\subsection*{Author Contributions}

All authors have participated in (a) conception and design, or analysis and interpretation of the data; (b) drafting the article or revising it critically for important intellectual content; and (c) approval of the final version.



\section*{Data availability statement}

The data that support the findings of this study are available from the corresponding author.

\end{document}